\documentclass[12pt]{article}
\usepackage{amsmath,amsthm,amscd,amsfonts}

\newcounter{num}
\newtheorem{theorem}{Theorem}[section]

\newtheorem*{Moore}{Moore's Conjecture}

\newtheorem*{Serre}{Serre's Theorem}

\newtheorem{definition}[theorem]{Definition}
\newtheorem{corollary}[theorem]{Corollary}
\newtheorem{proposition}[theorem]{Proposition}
\newtheorem{lemma}[theorem]{Lemma}

\begin{document}

\newcommand{\be}        {\begin{eqnarray}}
\newcommand{\ee}        {\end{eqnarray}}
\newcommand{\pl}{\partial}
\newcommand{\sbs}{\subset}
\newcommand{\vr}{\varphi}
\newcommand{\lm}{\lambda}
\newcommand{\eps}{\varepsilon}
\newcommand{\nb}{\nabla}
\newcommand{\wt}{\widetilde}

\newcommand{\cV}      {{\cal V}}
\newcommand{\cE}      {{\cal E}}
\newcommand{\cD}      {{\cal D}}
\newcommand{\cK}      {{\cal K}}
\newcommand{\cB}      {{\cal B}}
\newcommand{\cI}      {{\cal I}}
\newcommand{\cM}      {{\cal M}}
\newcommand{\cA}      {{\cal A}}
\newcommand{\cR}      {{\cal R}}
\newcommand{\cT}      {{\cal T}}
\newcommand{\cP}      {{\cal P}}
\newcommand{\cC}      {{\cal C}}
\newcommand{\cQ}      {{\cal Q}}
\newcommand{\cG}      {{\cal G}}
\newcommand{\cW}      {{\cal W}}
\newcommand{\cL}      {{\cal L}}
\newcommand{\cF}      {{\cal F}}
\newcommand{\cS}      {{\cal S}}
\newcommand{\cH}      {{\cal H}}
\newcommand{\cN}      {{\cal N}}

\newcommand{\AAA}       {{\Bbb A}}     
\newcommand{\aaa}       {{\LBbb A}}
\newcommand{\BB}        {{\Bbb B}}
\newcommand{\bb}        {{\LBbb B}}
\newcommand{\CC}        {\mathbb{C}}
\newcommand{\cc}        {{\LBbb C}}
\newcommand{\DD}        {{\Bbb D}}
\newcommand{\dd}        {{\LBbb D}}
\newcommand{\EEE}       {{\Bbb E}}
\newcommand{\eee}       {{\LBbb E}}
\newcommand{\FF}        {{\Bbb F}}
\newcommand{\ff}        {{\LBbb F}}
\newcommand{\GGG}       {{\Bbb G}}
\newcommand{\ggg}       {{\LBbb G}}
\newcommand{\HH}        {{\Bbb H}}
\newcommand{\hh}        {{\LBbb H}}
\newcommand{\II}        {{\Bbb I}}
\newcommand{\ii}        {{\LBbb I}}
\newcommand{\JJ}        {{\Bbb J}}
\newcommand{\jj}        {{\LBbb J}}
\newcommand{\KK}        {{\Bbb K}}
\newcommand{\kk}        {{\LBbb k}}
\newcommand{\LLL}       {{\Bbb L}}
\newcommand{\lll}       {{\LBbb L}}
\newcommand{\MM}        {{\Bbb M}}
\newcommand{\mm}        {{\LBbb M}}
\newcommand{\NN}        {\mathbb{N}}
\newcommand{\nn}        {{\LBbb N}}
\newcommand{\OO}        {{\Bbb O}}
\newcommand{\oo}        {{\LBbb O}}
\newcommand{\PP}        {{\Bbb P}}
\newcommand{\pp}        {{\LBbb P}}
\newcommand{\QQ}        {{\Bbb Q}}
\newcommand{\qq}        {{\LBbb Q}}
\newcommand{\RR}       { \mathbb{R}}
\newcommand{\rr}        {{\LBbb R}}
\newcommand{\SSS}       {{\Bbb S}}
\newcommand{\sss}       {{\LBbb S}}
\newcommand{\TTT}       {{\Bbb T}}
\newcommand{\ttt}       {{\LBbb t}}
\newcommand{\UU}        {{\Bbb U}}
\newcommand{\VV}        {{\Bbb V}}
\newcommand{\vv}        {{\LBbb V}}
\newcommand{\WW}        {{\Bbb W}}
\newcommand{\ww}        {{\LBbb W}}
\newcommand{\XX}        {{\Bbb X}}
\newcommand{\xx}        {{\LBbb X}}
\newcommand{\YY}        {{\Bbb Y}}
\newcommand{\yy}        {{\LBbb Y}}
\newcommand{\ZZ}        {\mathbb{Z}}
\newcommand{\zz}        {{\LBbb Z}}

\title{Profinite groups, profinite completions and a conjecture of
Moore}
\author{{Eli Aljadeff}\thanks{Supported by the Fund for the Promotion of Research at the Technion.}\\
Department of Mathematics\\
Technion -- Israel Institute of Technology\\
32000 Haifa, Israel\\
{\small e-mail: aljadeff@math.technion.ac.il; aljadeff@tx.technion.ac.il}}



\maketitle

\vspace{1cm}

\begin{abstract}
Let $R$ be any ring (with $1)$, $\Gamma $ a group and $R\Gamma $ the corresponding group ring. Let
$H$ be a subgroup of $\Gamma $ of finite index. Let $M$ be an $R\Gamma -$module, whose restriction
to $RH$ is projective.

Moore's conjecture \cite{C}: Assume for every nontrivial element $x$ in $\Gamma $, at least one of
the following two conditions holds:

$M1)$ \ $\langle x\rangle \cap H \neq \{e\}$ (in particular this holds if $\Gamma $ is torsion
free)

$M2)$ \ $ord(x)$ is finite and invertible in $R$. \\
Then $M$ is projective as an $R\Gamma$-module.

More generally, the conjecture has been formulated for crossed products $R\ast \Gamma $ and
even for strongly graded rings $R(\Gamma )$. 
We prove the conjecture for new families of groups, in particular for groups whose profinite
completion is torsion free.
 \medskip

The conjecture can be formulated 
for profinite modules $M$ over complete groups rings $[[R\Gamma ]]$ where $R$ is a profinite ring
and $\Gamma $ a profinite group. We prove the conjecture for arbitrary profinite groups. This
implies Serre's theorem on cohomological dimension of profinite groups.

\end{abstract}

\vspace{1cm}

\section{Introduction.}

Let $\Gamma $ be a group and $H$ a subgroup of finite index. Let ${\ZZ}\Gamma $ be the integral
group ring and let $M$ be a ${\ZZ}\Gamma -$module. It is well known (and easy to see) that if $M$
is projective over ${\ZZ}\Gamma $ then it is projective as a module over the subring ${\ZZ}H$ (in
fact this is true for arbitrary subgroups $H$ of $\Gamma ).$ Is the converse true? It is not
difficult to construct examples which show that this is false. For instance let $\Gamma $ be a
nontrivial finite group, $H=\{e\}$ and $M\cong {\ZZ}$ with the trivial $\Gamma -$action. Clearly
$M$ is ${\ZZ}={\ZZ}H-$free but not projective over ${\ZZ}\Gamma .$ More generally, one can
construct examples (that show that the converse is false) in which $\Gamma $ contains a nontrivial
element $x$ of finite order such that $<x>\cap \;H=\{e\}.$

Let $R$ be an arbitrary ring with unit element $1$ and let $\Gamma $ be any group. Let $R(\Gamma )$
be a strongly graded ring and $M$ a module over $ R(\Gamma ).$ Recall that if $M$ is projective
over $R(\Gamma )$ then $M$ is projective over the subring $R(U)$ where $U$ is any subgroup of
$\Gamma .$ Let $H$ be a finite index subgroup of $\Gamma $.


\begin{Moore}[see {\cite[Conjecture~1.1]{ACGK}}]\label{Moore}
  Assume for every nontrivial element $x$ in $\Gamma $, at least one of the following two
conditions holds:

$M1)$\ $\langle x\rangle \cap H \neq \{e\}$

$M2)$ \ $ord(x)$ is finite and invertible in $R.$\\
Then every $R(\Gamma )-$module $M$ which is projective over $R(H),$ is projective also over
$R(\Gamma ).$
\end{Moore}

We will say that Moore's conjecture holds for a group $\Gamma $ if the conjecture above holds for
an arbitrary ring $R$ and an arbitrary subgroup of finite index $H$ (see definitions in section 1).

One of the main advantages of strongly graded rings (over group rings) is that they ``allow
induction''. Recall that if $H$ is normal in $\Gamma $ then the strongly graded ring $R(\Gamma )$
can be expressed as a strongly graded ring of $R(H)$ with $\Gamma /H.$ (A similar result applies
for crossed products but obviously not for group rings). This ``flexibility'' will be used in our
proofs.

\medskip
\noindent
 {\it Remarks.}
\begin{enumerate}

\item The conjecture was formulated for group rings in \cite{C}. In \cite[Prop. 8]{C}  it is shown
that if Moore's conjecture holds for a class of groups $C$, then it holds for the class $LC$
(locally $C$) restricted to finitely generated modules. The main result in \cite {C} implies that
the conjecture holds for finite groups and hence the conjecture holds for locally finite groups
restricted to finitely generated modules.

\item If $R={\ZZ}$ the ring of integers and $\Gamma $ is torsion free, the conjecture says that
every ${\ZZ}\Gamma $-module $M$ which is ${\ZZ}H$-projective is also ${\ZZ}\Gamma $-projective.

\item Moore's conjecture is a far reaching generalization of Serre's theorem on cohomological
dimension of groups \cite {Swan}.

\begin{Serre}\label{Serre}
Let $\Gamma $ be a group and $H$ a subgroup of finite index. Assume $H$ has finite cohomological
dimension (that is $\Gamma $ has virtual finite cohomological dimension). If $\Gamma $ is torsion
free then it has finite cohomological dimension. Moreover, $cd(\Gamma )=cd(H).$
\end{Serre}

Indeed, Moore's conjecture implies Serre's theorem as follows: let $cd(H)=n$. If $P.$ $\rightarrow
{\ZZ}\rightarrow 0$ is a projective resolution of ${\ZZ}$ over ${\ZZ}\Gamma ,$ it is projective
also over $H.$ It follows that the $n$-th syzygy $Y_{n}$ of the resolution is a ${\ZZ}\Gamma $-
module whose restriction to $H$ is projective. Moore's conjecture says that $ Y_{n}$ is projective
over ${\ZZ}\Gamma $ and so $cd(\Gamma )\leq n.$

\item Moore's conjecture holds for groups which belong to a certain class of groups $H_{1}F.$
 The class $H_{1}F$ contains and is strictly larger than the class of groups of virtual
finite cohomological dimension. If $\Gamma \in H_{1}F$ and torsion free, then $\Gamma $ has finite
cohomological dimension. We refer the reader to (\cite[2.1]{K}) for the precise definition of the
class $ H_{1}F.$

\item Kropholler has constructed also a much larger class of groups, denoted by $HF.$ The class
$HF$ is extension closed, subgroup closed and closed under directed unions. In particular it
contains
\begin{enumerate}
\item  $H_{1}F$;

\item  every finitely generated soluble group;

\item  every countable linear group.
\end{enumerate}

\end{enumerate}

For this class one has


\begin{theorem}[\cite{ACGK}]\label{Th1.1}
Let $\Gamma \in HF$ and let $R(\Gamma )$ be a strongly graded ring over $R.$ Then Moore's
conjecture holds whenever the module $M$ is finitely generated.

\end{theorem}

Our first task is to ``reduce'' the problem to finitely generated groups (see also \cite[Prop.
8]{C} ). More precisely

\begin{theorem}\label{Th1.2}
If Moore's conjecture holds for every finitely generated subgroup of a group $G$ then it holds for
$G.$

\end{theorem}

The proof is based on a result of Benson and Goodearl (see section 2).

\begin{corollary}[{\cite[Cor.\,5.1]{C}}]\label{C1.3}
Moore's conjecture holds for any abelian group.
\end{corollary}

\begin{proof}
Every finitely generated abelian group belongs to $H_{1}F.$
\end{proof}

The main idea in this paper is to analyze the group $\Gamma $ from its $``top" $ rather than its
$``bottom''$ as in the construction of the classes $ H_{1}F$ and $HF.$ By $``top"$ we mean the
finite quotients of $\Gamma $ and by ``$bottom''$ we mean the subgroups of $\Gamma $ that appear as
``stabilizers'' in Kropholler's construction.

Let $P(\Gamma )$ be the collection of all finite index, normal subgroups of $ \Gamma .$ Let $\Omega
_{\Gamma }$ be a subset of $P(\Gamma )$ filtered from below. Assume further that $\Omega _{\Gamma
}$ is cofinal in $P(\Gamma ).$ Denote by $\widehat{\Gamma }={\lim }\Gamma /N$ the profinite
completion of $\Gamma $ with respect to $\Omega _{\Gamma }$ and let $\phi :\Gamma \rightarrow
\widehat{\Gamma }$ be the canonical map induced by the natural projections $\Gamma \rightarrow
\Gamma /N,\;N\in \Omega _{\Gamma }.$ Our main result (for abstract groups) is Theorem~\ref{Th3.1}.
Its formulation requires some terminology which is introduced in section~2. Its main corollary is

\begin{theorem}\label{Th1.4}
Let $\Gamma $ be a group and $\widehat{\Gamma }$ its profinite completion as above. If
$\widehat{\Gamma }$ is torsion free then Moore's conjecture holds for $\Gamma .$ More generally: If
any element $z$ of prime order (say $p)$ in $\widehat{\Gamma }$ is conjugate to an element $\phi
(x)$ where $x$ is an element (in $\Gamma $) of order $p,$ then Moore's conjecture holds for $\Gamma
.$
\end{theorem}

The condition in Theorem~\ref{Th1.4} is known to hold for large families of groups. For instance it
holds for families such as (see \cite {KW})
\begin{enumerate}
\item soluble minimax groups (see \cite{R} for the definition);

\item torsion free, finitely generated abelian by nilpotent groups.
\end{enumerate}

\noindent
 {\it Remark.}~~The finitely generated condition in 2) is important. Indeed, in \cite {KW}
an example is given of a residually finite, torsion free abelian by nilpotent group with torsion in
its profinite completion. Nevertheless, Moore's conjecture holds for such a group by
Theorem~\ref{Th1.2}.

It is not difficult to construct examples of finitely generated torsion free\linebreak
  metabelian
groups with infinite cohomological dimension whose profinite completion is torsion free (e.g.
$\Gamma ={\ZZ}^{\infty }\wr {\ZZ}$ $ =\{((x_{i})_{i\in {\ZZ}},\sigma ^{j})$ where $\sigma
((x_{i})_{i\in {\ZZ}})=(x_{i-1})_{i\in {\ZZ}}$ (right shift)). Note that $\Gamma $ does not belong
to $H_{1}F.$

The finitely generated groups considered so far belong to $HF.$ It is known that the Thompson group
$$
T=\langle x_{o},x_{1},x_{2},\dots,x_{n},\ldots:x_{n}^{x_{i}}=x_{n+1} \ \ {\rm for \ every} \ \
i<n\rangle
$$
 is not in $HF$, (see \cite{BrG, K}). On the other hand it is known that
\begin{enumerate}
\item $rad(T)=\cap \{H:[\Gamma:H]<\infty\}=T^{^{\prime }}$, the commutator subgroup of $T.$

\item $T_{ab}=T/T^{^{\prime }}\;$ is free abelian of rank 2.
\end{enumerate}

\begin{corollary}\label{C1.5}
Moore's conjecture holds for the Thompson group.
\end{corollary}

\begin{proof}
$\widehat{T}\cong \widehat{{\ZZ}^{2}}$ is torsion free.
\end{proof}

It was an open question whether the profinite completion $\widehat{\Gamma }$ of a torsion free,
residually finite group $\Gamma $ is necessarily torsion free. A first counterexample to this
question was given by Evans in \cite{E}.  Later, Lubotzky in \cite{Lu} gave an example of a torsion
free residually finite group whose profinite completion contains copies of any finite group!

Theorem~\ref{Th1.4} may be applied also to profinite groups (viewed as abstract group). Recently,
Nikolov and Segal announced the following important result.

\begin{theorem}[\cite{NS}]\label{Th1.6}
Let $\Gamma $ be a topologically finitely generated profinite group (this means that $\Gamma $ has
a f.g. subgroup which is dense in $ \Gamma $). Then every subgroup of finite index in $\Gamma $ is
open.
\end{theorem}

It follows that $\Gamma $ is naturally isomorphic to its profinite completion. In particular, the
condition in Theorem~\ref{Th1.4} is satisfied by $\Gamma $ and so we have the following

\begin{corollary}\label{C1.7}
Let $\Gamma $ be a profinite group, topologically finitely generated. Then Moore's conjecture holds
for $\Gamma .$
\end{corollary}

Using Theorem~\ref{Th1.2} we obtain Moore's conjecture for arbitrary profinite groups considered as
abstract groups. We record this in

\begin{corollary}\label{C1.8}
Let $\Gamma $ be a profinite group. Then Moore's conjecture holds for $\Gamma $ (as an abstract
group).
\end{corollary}

The interest in Corollaries~\ref{C1.7} and \ref{C1.8} is limited. It is of course desired to obtain
similar results in the category of profinite modules. For this one needs to adapt Chouinard's
results in \cite {C}, to profinite rings and profinite modules. Let $R$ be a profinite ring and
$\Gamma $ a profinite group. Denote by $[[R\Gamma ]]$ the complete group ring (see \cite{RZ}
section 5.3). Our main result in this context (and perhaps in the entire paper) is

\begin{theorem}\label{Th1.9}
Let $[[R\Gamma ]]$ be a complete group ring and $M$ an $[[R\Gamma ]]$-profinite module. Let $H$ be
an open subgroup of $\Gamma .$ Assume Moore's condition $(M1$ or $M2)$ holds$.$ Then if $M$ is
projective over $[[RH]]$ then $M$ is projective over $[[R\Gamma ]].$
\end{theorem}

This implies

\begin{corollary}\label{C1.10}
Let $[[R\Gamma ]]$ be a complete group ring and $H$ an open subgroup of $\Gamma .$ Assume Moore's
condition ($M1$ or $M2)$ holds. Then for any $[[R\Gamma ]]$-profinite module $M$ we have $proj.dim
_{[[R\Gamma ]]}(M)=proj.dim _{[[RH]]}(M).$
\end{corollary}

As a direct consequence we get Serre's theorem on cohomological dimensions of profinite groups (See
\cite{Se},\cite{Ha}.)

\begin{theorem}\label{Th1.11}
Let $\Gamma $ is be a profinite group and $H$ an open subgroup. If $\Gamma $ has no elements of
order $p$ then $cd_{p}(\Gamma )=cd_{p}(H).$
\end{theorem}

This paper is organized as follows: In section~2 we set most of the terminology and notation needed
in the paper. The section contains also two reductions one of which is the reduction to finitely
generated subgroups mentioned above. The other reduction allows us to replace the subgroup of
finite index $H$ in Moore's conjecture with its core in $\Gamma ,$ namely the intersection of all
its conjugates in $\Gamma .$ Next, we continue with a brief discussion on the necessity of the
condition in Moore's conjecture (basically the condition is necessary only if $H$ is normal in
$\Gamma )$. We close the section by recalling some basic facts on profinite topologies.

The main results are in sections 3 and 4. In section 3 we prove Theorem~\ref{Th1.4}. As mentioned
above it follows from Theorem~\ref{Th3.1} which contains the main construction in the paper. The
last section, section 4, contains the proof of Moore's conjecture for complete group rings and
profinite modules.


\section{Preliminaries, Terminology and Reductions.}

\setcounter{equation}{0}
\setcounter{theorem}{0}

It is convenient to use the following terminology: if condition $(M1$ or $M2)$ holds for the group
$\Gamma $, the subgroup of finite index $H$ and the coefficient ring $R,$ we will say that Moore's
condition holds for the triple $(\Gamma ,H,R).$ If $M1$ holds for the group $%
\Gamma $ and the subgroup $H$ we will say that Moore's condition holds for the pair $(\Gamma ,H).$
Note that $M1$ holds for $(\Gamma ,H)$ if and only if $(M1$ or $M2)$ holds for $(\Gamma ,H,R),$
where $R$ is any ring. Next we'll say that Moore's conjecture holds for the triple $(\Gamma ,H,R)$
if condition $(M1$ or $M2)$ implies that any module $M$ over $R(\Gamma )$ which is projective over
$R(H),$ is projective over $R(\Gamma ).$ We'll say that Moore's conjecture holds for $(\Gamma ,H)$
if condition $M1$ implies the same conclusion for $(\Gamma ,H,R)$ where $R$ is arbitrary. Finally,
we'll say that Moore's conjecture holds for the group $\Gamma $ if the conjecture holds for
$(\Gamma ,H)$ , for any subgroup $H$ of finite index in $\Gamma .$

We turn now to the proof of Theorem~\ref{Th1.2} starting with the following lemma which is a
particular case of \cite[p.125 exercise 17]{CE}.

\begin{lemma}\label{L2.1}
Let $G$ be any group and $R(G)$ a strongly graded ring over $R.$ Let $M$ be a module over $R(G).$
If $M$ is flat over any subring of the form $R(\Gamma )$ where $\Gamma $ is a finitely generated
subgroup of $G$ then $M$ is flat over $R(G).$
\end{lemma}

As mentioned in the introduction, the proof of Theorem~\ref{Th1.2} is based on a theorem of Benson
and Goodearl. For the reader's convenience we recall it here:

\begin{theorem}[\cite{BeG}]\label{Th2.2}
Let $R(G)$ be a strongly graded ring over $R$ and let $H$ be a subgroup of finite index of $G.$ Let
$M$ be a flat module over $R(G)$ which is projective over $R(H).$ Then $M$ is projective over
$R(G).$
\end{theorem}

\begin{proof}[Proof of Theorem~\ref{Th1.2}.]
Let $H$ be a subgroup of finite index of $G$ and let $R$ be any ring. Assume Moore's condition
holds for the triple $(G,H,R)$ and let $M$ be a module over $R(G),$ projective over $R(H).$ We need
to show $M$ is projective over $R(G).$ Applying Theorem~\ref{Th2.2} it is sufficient to prove that
$M$ is flat over $R(G).$ By Lemma~\ref{L2.1} it is sufficient to show that $M$ is flat over
$R(\Gamma )$ for every finitely generated subgroup $\Gamma $ of $ \;G.$ Let $\Gamma $ be such a
group and let $H_{\Gamma }=H\cap \Gamma .$ Clearly,$\;M$ is projective over $R(H_{\Gamma }).$
Furthermore, since Moore's condition holds for the triple $(G,H,R)$ it holds also for the triple
$(\Gamma ,H_{\Gamma },R)$ and hence $M$ is projective over $R(\Gamma ).$ This implies that $M$ is
flat over $R(\Gamma )$ and the result follows.
\end{proof}

Two important ingredients in the proof of the main theorem (Theorem~\ref{Th3.1}) are Chouinard's
theorem and Maschke's theorem for strongly graded rings. Since we will be using them repeatedly we
recall them here starting with Chouinard's

\begin{definition}\label{D2.3}
Let $\Gamma $ be a finite group. We say that an $ R(\Gamma )$-module $M$ is weakly projective if
there is an $f\in Hom_{R}(M,M) $ with $tr_{\Gamma }(f)=\sum_{\sigma \in \Gamma }\sigma
(f)=id_{M}$\;(here $\Gamma $ acts on $Hom_{R}(M,M)$ diagonally).
\end{definition}

\begin{theorem}[\cite{AG}]\label{Th2.4}
Let $R$ be an arbitrary ring with identity and let $\Gamma $ be a finite group. Let $R(\Gamma )$ be
a strongly graded ring. If $M$ is any module over $R(\Gamma )$ then it is weakly projective
(projective) if and only if it is weakly projective (projective) over all subrings $R(P)$ where $
P$ is an elementary abelian subgroup of $\Gamma.$ In fact it is sufficient to assume that $M$ is
projective over $R(P)$ where $P$ runs over representatives of all conjugacy classes of maximal
elementary abelian subgroups of $\Gamma .$
\end{theorem}

\begin{theorem}[Maschke, {\cite[Chapter 1,\,section 4]{Pass}}]\label{Th2.5}
Let $R$ be an arbitrary ring with identity and let $\Gamma $ be a finite group whose order is
invertible in $R$. Let $R(\Gamma )$ be a strongly graded ring. Then any $R(\Gamma )-$ module is
weakly projective. In particular, any module $M$ over $R(\Gamma )$ is projective if and only if it
is projective over $R.$
\end{theorem}

\begin{quote}
\begin{itemize}
\item
Note that Chouinard's and Maschke's theorems imply Moore's conjecture for $(\Gamma ,H,R)$ whenever
the group $\Gamma $ is finite.
\end{itemize}
\end{quote}

\medskip

In the proof of Theorem~\ref{Th3.1} it is convenient to replace $H$ by a normal subgroup $H_{0}.$
Let $H$ be a subgroup of finite index in $\Gamma .$ Let $ H_{0}$ be the $core$ of $H$ in $\Gamma $
that is $H_{0}=\cap H^{g}.$

\begin{lemma}\label{L2.6}
Moore's condition holds for the triple $(\Gamma ,H,R)$ if and only if it holds for the triple
$(\Gamma ,H_{0},R).$ Furthermore, if Moore's conjecture holds for the triple $(\Gamma ,H_{0},R)$
then it holds for the triple $ (\Gamma ,H,R).$
\end{lemma}

\begin{proof}
If $z$ is of order $p$ and not in $H_{0}$ then $z$ is not in one of the conjugates $H^{g}$ of $H.$
Then $z^{g^{-1}}$ which is of order $p,$ is not in $H.$ The second statement follows from the fact
that a projective module over $R(H)$ is projective over $R(H_{0}).$
\end{proof}

It is natural to ask whether Moore's condition for a pair $(\Gamma ,H)$ is necessary in Moore's
conjecture. More precisely we ask the following questions:
\begin{enumerate}
\item
Assume Moore's condition does not hold for $(\Gamma ,H).$ Is there a ring $R,$ a strongly graded
ring $R(\Gamma )$ and a module $M$ over $ R(\Gamma )$ which is projective over $R(H)$ but not
projective over $ R(\Gamma )?$
\item
 Let $R$ be given and assume Moore's condition does not hold for $(\Gamma ,H,R).$ Is there a
module $M$ over a strongly graded ring $R(\Gamma )$ which is projective over $R(H)$ but not over
$R(\Gamma )$?
\end{enumerate}

\begin{proposition}\label{P2.7}
If the group $H$ is normal in $\Gamma $ then Moore's condition is necessary in the following
(strong) sense: If Moore's condition does not hold for $ (\Gamma ,H,R)$ then there is a module $M$
over $R\Gamma $ (the group ring) which \ is projective over $RH$ but not over $R\Gamma .$
\end{proposition}

\noindent
 {\it Remark.}~~We cannot expect the proposition to hold for every strongly graded ring
($e.g.$ the Gauss integers ${\ZZ}[i]$ may be represented as a twisted group ring of the group of
two elements over ${\ZZ}.$ Moore's condition does not hold for $(C_{2},\{e\},{\ZZ})$ but every
module over ${\ZZ}[i]$ which is projective over $ {\ZZ}$ is projective over ${\ZZ}[i]).$
\medskip

\begin{proof}[Proof of Proposition~\ref{P2.7}.]
Let $\sigma $ be an element in $\Gamma \backslash  H$ of order $p$ where $p$ is not invertible in
$R.$ Denote by $U$ the cyclic group it generates. Consider the left $R\Gamma $-module $M=R\Gamma
\otimes _{{\ZZ}U}{\ZZ}$ (with the obvious left $R\Gamma $-structure). The module $M$ is free over
$RH$ with a basis consisting a set of representatives for the cosets of the subgroup $\langle
H,U\rangle$ in $\Gamma.$ We claim that the map $\pi :R\Gamma \rightarrow M,\;\pi (g)=g\otimes 1$
does not split over $R\Gamma.$ Indeed, if $j$ is a splitting over $R\Gamma $, let
$$
j(1\otimes 1)=\sum\limits_{i=0}^{p-1}\alpha _{i}\sigma
^{i}+\sum\limits_{i=0}^{p-1}\sum\limits_{s_j\in T_{[\Gamma :U]}\backslash \{1\}}^{}
 \alpha _{ji}s_{j}\sigma ^{i}
$$
where $T_{[\Gamma :U]}$ is a set of representatives for the cosets of $U$ in $\Gamma $ and such
that $1\in T_{[\Gamma :U]}.$ By the splitting condition we get that
$$
\sum\limits_{i=0}^{p-1}\alpha _{i}=1.
 $$

Furthermore by the $R\Gamma $-linearity of the map $j$ we get that $\alpha _{0}=\dots=\alpha
_{p-1}$ and so $p$ is invertible in $R.$ Contradiction.
 \end{proof}

The fact that the group $H$ is normal in $\Gamma $ is essential in Proposition~\ref{P2.7}. Next we
give an example of a finite group $\Gamma $ and a subgroup $H$ such that $\Gamma\backslash H$
contains an element of prime order, but for any ring $R$, any strongly graded ring $R(\Gamma )$\
and any module $M$ over $ R(\Gamma ),$ if $M$ is projective over $R(H),$ then it is projective also
over $R(\Gamma ).$

\noindent
 {\it Example.}~~Let $\Gamma =\langle\sigma ,\tau :\sigma ^{9}=\tau ^{2}=1,\tau \sigma \tau
=\sigma ^{-1}\rangle\;$the dihedral group of order $18.$ Let $H=\langle\sigma ^{3},\tau \rangle$ .
The set $\Gamma \;\backslash H$ contains elements of order $2$ and hence Moore's condition does not
hold for $(\Gamma ,H).$ On the other hand every elementary abelian subgroup of $\Gamma $ is cyclic
and conjugate to a subgroup of $H.$ The result follows from Chouinard's theorem.

\medskip

We close this section by recalling some basic concepts and fixing notation concerning profinite
topologies and profinite completions of groups (we refer the reader to \cite {RZ}). A non-empty
collection $\Sigma $ of normal subgroups of finite index of a group $\Gamma $ is filtered from
below if for any $N_{1},N_{2}\in \Sigma $ there exists $N\in \Sigma $ such that $N\subseteq
N_{1}\cap N_{2}.$ Then $\Gamma $ turns into a topological group by considering $\Sigma $ as a
fundamental system of neighborhoods of the identity element $1$ of $\Gamma.$ We denote by
$K_{\Sigma}(\Gamma)= {\lim }\Gamma /N $ the profinite completion with respect to that topology and
by $\phi :\Gamma \rightarrow $ $K_{\Sigma }(\Gamma )$ the canonical map induced by the natural
projections $\Gamma \rightarrow \Gamma /N,\;N\in \Sigma.$

Given any collection $\Sigma $ as above let $I_{\Sigma }$ be the index set that correspond to
$\Sigma ,$ that is for every $N\in \Sigma $ we have $ i_{N}\in I_{\Sigma }.$ Obviously, using the
ordering $i_{N_{1}}>i_{N_{2}}$ if and only if $N_{1}<N_{2},$ the set $I_{\Sigma }$ is partially
ordered and directed. For $k<j$ in $I_{\Sigma }$, we denote by $\phi _{jk}:\Gamma /N_{j}\rightarrow
\Gamma /N_{k}$ the natural projection.

Finally, recall that if $H$ is a subgroup of $\Gamma $ of finite index and $ \Sigma $ is a
collection as above with the additional condition that all elements $N\in \Sigma $ are contained in
$H,$ then there is a natural inclusion $K_{\Sigma }(H)\hookrightarrow K_{\Sigma }(\Gamma ).$


\section{Abstract groups.}

\setcounter{equation}{0}
\setcounter{theorem}{0}

\begin{theorem}\label{Th3.1}
Let $\Gamma $ be any group and $H$ a subgroup of finite index. Let $\Sigma $ be a collection of
subgroups of $H,$ filtered from below, normal and of finite index in $\Gamma .$ Let $K_{\Sigma
}(\Gamma )$ and $K_{\Sigma }(H)$ be the profinite completions of $\Gamma $ and $H$ with respect to
$\Sigma.$ Assume any element of prime order (say $p)$ in $K_{\Sigma }(\Gamma )\;\backslash
\;K_{\Sigma }(H)$ is conjugate to an element $\phi (x)$ where $ x$ is an element in $\Gamma $ of
order $p.$ Then Moore's conjecture holds for $(\Gamma ,H).$ In particular if $K_{\Sigma }(\Gamma )$
is torsion free then Moore's conjecture holds for $(\Gamma ,H).$
\end{theorem}

\begin{proof}
Let $H_{0}$ be the core of $H$ in $\Gamma .$ By Lemma~\ref{L2.6} it is sufficient to show the
conjecture for $(\Gamma ,H_{0}).$ Note that since the groups in $\Sigma $ are normal in $\Gamma $,
they are contained in $H_{0}.$ For every $j\in I_{\Sigma },$ we denote by $\phi _{j}:\Gamma
/N_{j}\rightarrow \Gamma /H_{0}$ the natural projection.

For every $j\in I_{\Sigma },$ let $Y_{j}$ be the subset of $\prod_{i\in I_{\Sigma }}\Gamma /N_{i}$
defined by
$$
 Y_{j}=\{(x_{i})\in \prod_{i\in I_{\Sigma }}\Gamma /N_{i} : \phi _{jk}(x_{j})=x_{k} \ \ {\rm
 whenever} \ \ k<j\}.
$$

For any prime number $p$,  let $ Z_{j}^{(p)}\subset Y_{j}$ be the set
 $$
Z_{j}^{(p)}=\{(x_{i})\in Y_{j}:ord(x_{j})=p \ \ {\rm in} \ \ \Gamma /N_{j} \ \ {\rm and} \ \
ord(\phi _{j}(x_{j}))=p \ \  {\rm in} \ \  \Gamma /H_{0}\}.
$$

Note that $Z_{j}^{(p)}\supset Z_{j^{^{\prime }}}^{(p)}$ for $j<j^{^{\prime }}.$

The main step of the proof is the first statement in the following lemma. Recall that by definition
Moore's conjecture holds for $(\Gamma ,H_{0})$ if and only if Moore's conjecture holds for $(\Gamma
,H_{0},R)$ with $R$ arbitrary.
\end{proof}

\begin{lemma}\label{L3.2}
If the conclusion in Moore's conjecture does not hold for $(\Gamma ,H_{0},R)$ (that is there is a
non-projective $R(\Gamma )$-module $M$ which is projective over $R(H_{0})),$ then there is a prime
number $p,$ not invertible in $R$, such that $Z_{j}^{(p)}$ is non-empty for every $j\in I_{\Sigma
}.$ Furthermore, $Z_{j}^{(p)}$ is closed in $\;\prod_{i\in I_{\Sigma }}\Gamma /N_{i}.$
\end{lemma}

Let us postpone the proof of the lemma and complete first the proof of the theorem. Assume the
theorem is false. This means that Moore's conjecture does not hold for a triple $(\Gamma
,H_{0},R),$ some $R.$ In particular the conclusion in Moore's conjecture does not hold for $(\Gamma
,H_{0},R)$. Let $ p$ be the prime number given by the lemma. From the condition $
Z_{j}^{(p)}\supset Z_{j^{^{\prime }}}^{(p)}$ for $j<j^{^{\prime }}$ it follows that the family
$\{Z_{j}^{(p)}\}_{j\in I_{\Sigma }}$ satisfies the finite intersection property. Moreover, since
$\prod_{i\in I_{\Sigma }}\Gamma /N_{i}$ is compact and the sets $Z_{j}^{(p)}\;^{\prime }s$ are
closed, there is an element $z\in \cap _{i\in I_{\Sigma }}Z_{i}^{(p)}.$ Clearly, $z$ is an element
of order $p$ in $K_{\Sigma }(\Gamma ).$ Furthermore, by the definition of the sets $Z_{j}^{(p)},$
we have $z\notin K_{\Sigma }(H_{0}).$ Now, it is easily checked that the condition (in the theorem)
on elements of $K_{\Sigma }(\Gamma )\;\backslash \;K_{\Sigma }(H)$ holds for all elements in
$K_{\Sigma }(\Gamma )\;\backslash \;K_{\Sigma }(H_{0})$ and hence there is an element $x$ in
$\Gamma $ of order $p$ such that $\phi (x)$ is conjugate to $z$ (in $K_{\Sigma }(\Gamma )).$ It is
clear that $x\notin H_{0} $ and we get a contradiction to Moore's condition $(M1$ or $M2)$ for $
(\Gamma ,H_{0},R).$
\medskip

\begin{proof}[Proof of Lemma~\ref{L3.2}.]
 Let $M$ be an $R(\Gamma )$-module,
projective over $ R(H_{0})$ but not projective over $R(\Gamma ).$ Write $R(\Gamma )=R(H_{0})(\Gamma
/H_{0}).$ By Chouinard's and Maschke's theorems there is a prime number $p$, not invertible in $R,$
and an elementary abelian $p-$group $E$ in $\Gamma /H_{0}$ such that $M$ is not projective over the
ring $ R(H_{0})(E).$ Let $T$ be the inverse image of $E$ in $\Gamma $ (with respect to the group
extension $1\rightarrow H_{0}\rightarrow \Gamma \rightarrow \Gamma /H_{0}\rightarrow 1).$ Now, for
any $N_{j}\in \Sigma $ write $ R(T)=R(N_{j})(T/N_{j}).$ Since $N_{j}\subset H_{0},$ $\ M$ is
projective over $R(N_{j}).$ Furthermore, $M$ is not projective over $
R(T)=R(H_{0})(E)=R(N_{j})(T/N_{j})$ and hence by Chouinard's theorem there is an elementary abelian
subgroup $E_{j}$ of $T/N_{j}$ such that $M$ is not projective over $R(N_{j})(E_{j}).$ Let $T_{j}$
be the inverse image of $ E_{j} $ in $T.$ The module $M$ is not projective over $R(T_{j})$ and
therefore $T_{j}$ is not contained in $H_{0}.$ We conclude that $T_{j}$ $/N_{j}=E_{j}$ is
elementary $p-$abelian with nontrivial image modulo $H_{0}. $ This proves the first statement of
the lemma. For the second statement in the lemma observe that the set $Z_{j}^{(p)}$ contains $
\Gamma /N_{i}$ for all $i$ but a finite subset of $I_{\Sigma }.$ This completes the proof of the
lemma and hence of Theorem~\ref{Th3.1}.
 \end{proof}

We close this section with the
\medskip

\begin{proof}[Proof of Theorem~\ref{Th1.4}.]
 Assume Moore's condition holds for $(\Gamma ,H,R)$ where $H$
is asubgroup of finite index of $\Gamma $ and $R$ is any ring. By Lemma 1.5 it is sufficient to
prove that Moore's conjecture holds for $(\Gamma ,H_{0},R).$ Let $\Omega _{\Gamma
}(H_{0})=\{H_{0}\cap N:N\in \Omega _{\Gamma }\}.$ Since $K_{\Omega _{\Gamma }(H_{0})}(\Gamma )\leq
K_{\Omega _{\Gamma }}(\Gamma )$ (in fact they are naturally isomorphic) the condition in
Theorem~\ref{Th3.1} (for elements in $K_{\Omega _{\Gamma }(H_{0})}(\Gamma )\;\backslash \;K_{\Omega
_{\Gamma }(H_{0})}(H_{0})$) holds and the result follows.
\end{proof}


\section{Profinite \ groups, \  complete \ group \ rings \ and  \ profinite \ modules.}

\setcounter{equation}{0}
\setcounter{theorem}{0}

Let $\Gamma $ be a profinite group and $H$ an open subgroup of $\Gamma .$ We denote by $T_{[\Gamma
:H]}$ a transversal for the right cosets of $H$ in $ \Gamma $ (we choose $1\in \Gamma $ as the
representative of the trivial coset)$.$ Let $M$ be a profinite module over the complete group ring
$ [[R\Gamma ]].$ We say that $M$ is $([[R\Gamma ]],[[RH]])-$ relative projective if given any
diagram of $[[R\Gamma ]]$-modules and $[[R\Gamma ]]-$ continuous maps $(\alpha ,\phi )$
$$
\begin{aligned}
  \ &                          &&   M              &                 &       \\
  \ &   \swarrow               &&   \downarrow\phi &                 &        \\
B \ &  \longrightarrow         && A                & \longrightarrow & \ \  0 \\[-0.5cm]
  \ &  \hspace{0.3cm} \alpha   &&                  &\\
\end{aligned}
$$

\noindent
 then if $\phi $ can be lifted to a continuous map $M\rightarrow B$ over $ [[RH]]$ then it
can be lifted also to a continuous map $ v:M\rightarrow B$ over $[[R\Gamma ]].$

Maschke's theorem in this context reads:

\begin{theorem}\label{Th4.1}
If $[\Gamma :H]$ (the index of $H$ in $\Gamma $) is invertible in $R$ then every profinite module
$M$ over $[[R\Gamma ]]$ is $([[R\Gamma ]],[[RH]])-$ relative projective.
\end{theorem}

\begin{proof}
With the above notation, if $\theta _{H}:M\rightarrow B$ is a continuous lifting of $\phi $ over
$[[RH]]$ then $ord(T_{[\Gamma :H]})^{-1}\sum_{\sigma \in T_{[\Gamma :H]}}\sigma \theta _{H}\sigma
^{-1}$ is a continuous lifting of $\phi $ over $[[R\Gamma ]].$
\end{proof}

\begin{lemma}\label{L4.2}
With the above notation, a module $M$ is $([[R\Gamma ]],[[RH]])-$ relative projective if and only
if there is an $[[RH]]-$ continuous map $ s:M\rightarrow M$ such that $tr_{H\rightarrow
G}(s)=\sum_{\sigma \in T_{[\Gamma :H]}}\sigma (s)=id_{M}$ (where $\sigma $ acts on $
Hom_{[[RH]]}(M,M)$ diagonally).
\end{lemma}

\begin{proof}
Consider the induced module $[[R\Gamma]]\hat{\otimes} _{[[RH]]}M$ of $M$ where $\hat{\otimes }$
denotes the complete tensor product (see \cite[section 5.5]{RZ}). Since $[[R\Gamma ]]$ is finitely
generated as an $[[RH]]-$ module, the complete tensor product coincides with the usual tensor
product $[[R\Gamma ]]\otimes _{\lbrack \lbrack RH]]}M.$

Assume $M$ is $([[R\Gamma ]],[[RH]])-$relative projective and consider the diagram of $[[R\Gamma
]]-$ modules
$$
\begin{aligned}
       \ &                          &&   M              &                 &       \\
       \ &   \swarrow               &&   \downarrow\ id &                 &        \\
[[R\Gamma]]\otimes_{[[RH]]}M \ &  \longrightarrow         &&   M              & \longrightarrow & \ \  0 \\[-0.5cm]
       \ &  \hspace{0.3cm} \pi      &&                  &\\
\end{aligned}
$$
where $\pi $ is given by $\sigma \otimes m\rightarrow \sigma m.$ The map $ \pi $ splits over
$[[RH]]$ $(j:m\rightarrow 1\otimes m$ is a splitting) and so it splits over $[[R\Gamma ]].$ Let
$v:M\rightarrow $\ $[[R\Gamma ]]\otimes _{\lbrack \lbrack RH]]}M$ be a splitting of $\pi $ over
$[[R\Gamma ]].$ Since the elements of $T_{[\Gamma :H]}$ form a basis of $[[R\Gamma ]]$ over $
[[RH]]$ we may write
$$
 v(m)=\sum_{\sigma \in T_{[\Gamma :H]}}\sigma \otimes f_\sigma(m)
$$
 where $f_{\sigma }:M\rightarrow M$ is a well defined map.

Now, one checks that the map $f_{\sigma }$ is $R-$ linear and moreover $ f_{\sigma }(hm)=
\linebreak
 \sigma ^{-1}h\sigma f_{\sigma }(m)$ for every $h\in H$ and $ m\in M.$ In particular
$f_{e}:M\rightarrow M$ is an $RH-$linear map. We show that $f_{e}$ is continuous (and hence
$[[RH]]-$ linear). Let $M_{0}\subset M$ be an open $[[RH]]-$submodule. We claim that
$f_{e}^{-1}(M_{0})$ contains an open $[[RH]]-$ submodule of $M.$ To see this, observe that
$[[R\Gamma ]]\otimes _{\lbrack \lbrack RH]]}M_{0}$ is open in $[[R\Gamma ]]\otimes _{\lbrack
\lbrack RH]]}M$ and by the continuity of $v$ we have that $v^{-1}([[R\Gamma ]]\otimes _{\lbrack
\lbrack RH]]}M_{0})=M_{1}$ is open $M.$ By the linear independence of the elements in $T_{[\Gamma
:H]}$ over $[[RH]]$ we have that $f_{e}(M_{1})\subset M_{0}$ and the claim is proved. Finally we
show that $tr_{H\rightarrow G}(f_{e})=id_{M}. $ By the $\Gamma $-linearity of $v$ we obtain that
$f_{\sigma }(m)=f_{e}(\sigma ^{-1}m)$ and from the equality $\pi v=id_{M}$ we obtain
$tr_{H\rightarrow G}(f_{e})=\sum_{\sigma \in T_{[\Gamma :H]}}\sigma f_{e}\sigma^{-1}=id_{M}$.

For the converse let $f:M\rightarrow M$ be a continuous $[[RH]]-$ map such that $tr_{H\rightarrow
\Gamma }(f)=id_{M}$. Let
$$
\begin{aligned}
  \ &                          &&   M              &                 &       \\
  \ &   \theta_H\swarrow               &&   \downarrow\phi &                 &        \\
B \ &  \longrightarrow         && A                & \longrightarrow & \ \  0 \\[-0.5cm]
  \ &  \hspace{0.3cm} \alpha   &&                  &\\
\end{aligned}
$$
be a diagram of $[[R\Gamma ]]$-modules where $\phi $ and $\alpha $ are $[[R\Gamma ]]$-continuous
maps and $\theta _{H}$ an $[[RH]]$-continuous map. The map $\theta _{\Gamma }=tr_{H\rightarrow
\Gamma }(\theta _{H}f)$ is an $[[R\Gamma ]]$-continuous map and it is a lifting of $\phi .$ The
lemma is proved.
\end{proof}

Assume now, $H$ is normal (and open) in $\Gamma .$ Denote by
$$
 Hom_{[[RH]]}(M,M)=End_{[[RH]]}(M)
$$
 the endomorphism ring of all continuous $ [[RH]]-linear$
homomorphisms $\phi :M\rightarrow M.$ The quotient group $ \Gamma /H$ acts on $End_{[[RH]]}(M)$ via
the diagonal action. Recall that a $ \Gamma /H$-module $U$ is weakly projective if and only if
there is an additive map $\phi :U\rightarrow U$ with $\sum_{\sigma \in \Gamma /H}\sigma \phi \sigma
^{-1}=id_{U}.$

\begin{proposition}\label{P4.3}
Let $\Gamma $ be a profinite group, $H$ a normal open subgroup of $\Gamma $, $R$ a profinite ring
and $M$ a profinite $[[R\Gamma ]]-$ module. Then the following conditions are equivalent:

\begin{enumerate}
\item $M$ is $([[R\Gamma ]],[[RH]])-$relative projective.

\item  There exists a continuous $[[RH]]-$map $f:M\rightarrow M$ with $ tr_{H\rightarrow \Gamma
}(f)=id_{M}.$

\item  $Hom_{[[RH]]}(M,M)$ is weakly projective $\Gamma /H-$module (with the diagonal action).

\item  $Hom_{[[RH]]}(M,M)$ is weakly projective $E-$module for every elementary abelian subgroup
$E$ of $\Gamma /H.$
\end{enumerate}
\end{proposition}

\begin{proof}
The equivalence of $(1)$ and $(2)$ is in Lemma~\ref{L4.2}.
\medskip

 \noindent
 \underline{$(2\rightarrow 3)$}~~Let $f:M\rightarrow M$ be a continuous $[[RH]]-$map with $ \sum_{\sigma \in
\Gamma /H}\sigma (f)=id_{M}.$  Applying the functor $ Hom_{[[RH]]}(-,M)$ to $f$ we obtain a
map\linebreak $f^{\ast }:Hom_{[[RH]]}(M,M)\rightarrow Hom_{[[RH]]}(M,M)$ of abelian groups such
that\linebreak $\sum_{\sigma \in \Gamma /H}\sigma f^{\ast }\sigma ^{-1}=id_{Hom_{[[RH]]}(M,M)}.$
 \medskip

 \noindent
\underline{$(3\rightarrow 2)$}~~Let $\phi :Hom_{[[RH]]}(M,M)\rightarrow Hom_{[[RH]]}(M,M) $ be a
map of abelian groups such that $\sum_{\sigma \in \Gamma /H}\sigma \phi
\sigma^{-1}=id_{Hom_{[[RH]]}(M,M)}.$ Then $\phi (id_{M}):M\rightarrow M $ is a continuous
$[[RH]]-$map with $\sum_{\sigma \in \Gamma /H}\sigma (\phi (id_{M}))=id_{M}.$
 \medskip

\noindent
 \underline{$(3\leftrightarrow 4)$}~~(see \cite[Theorem\,1]{AG}).
\end{proof}

\begin{theorem}\label{Th4.4}
 Let $H$ be a normal open subgroup
of a profinite group $\Gamma .$ Let $M$ be a profinite $[[R\Gamma ]]-$module. Then $M$ is
projective over $ [[R\Gamma ]]$ if and only if $M$ is projective over $[[RT]]$ for every $ T\leq
\Gamma $ which contains $H$ and $T/H$ is elementary abelian.
\end{theorem}

\begin{proof}
Observe that an $[[R\Gamma ]]-$module $M$ is projective if and only if it is projective over
$[[RH]]$ and $([[R\Gamma ]],[[RH]])-$relative projective. The theorem now follows from
Proposition~\ref{P4.3} implications $ (1\rightarrow 3)$ for $(T,H),\;(4\rightarrow
3),\;(3\rightarrow 1)\;$for $ (\Gamma ,H).$
\end{proof}

We can prove now our main theorem for profinite groups (Theorem~\ref{Th1.9}).

\begin{proof}
By taking the core of $H$ in $\Gamma $ we can assume $H$ is normal in $\Gamma .$ Let $\Omega
_{\Gamma }$ denote the family of all normal, open subgroups of $\Gamma $ and let $\Omega _{\Gamma
}(H)=\{H\cap N:N\in \Omega _{\Gamma }\}.$ Clearly, the family $\Omega _{\Gamma }(H)$ is cofinal in
$ \Omega _{\Gamma }$ and hence the completions of $H$ and $\Gamma $ with respect to $\Omega
_{\Gamma }(H)$ are naturally isomorphic to $H$ and $ \Gamma $ respectively. In particular the
condition in Theorem~\ref{Th3.1} is satisfied in an obvious way. The proof of the theorem is
completed by following the steps in the proof of Theorem~\ref{Th3.1}. The only place here which is
slightly different is in the proof of Lemma~\ref{L3.2}. Here instead of using strongly graded
rings, we apply Theorem~\ref{Th4.4}. Details are left to the reader.
\end{proof}

\noindent
 {\it Final Remark.}~~The original proof of Serre of Theorem~\ref{Th1.11} uses Serre's result on
product of Bockstein operators \cite {Se}. It is not surprising that Chouinard's theorem (which is
the main ingredient in our proof) is based also on Serre's result on Bockstein operators.

\bigskip

\noindent
 {\bf Acknowledgment.}~~In the search for examples of (abstract) groups which satisfy the
condition in Theorem~\ref{Th1.4}, I consulted several people. I would like to thank M. Sageev, M.
Sapir and in particular A. Lubotzky for their contribution in this matter. I would like to thank R.
Holzman for useful conversations I had with him. Special thanks to my colleague J. Sonn for the
infinite number of helpful conversations we had.

\end{document}